\documentclass{article}
\usepackage{latexsym}
\usepackage{amssymb}
\usepackage[latin1]{inputenc}
\newcommand{\R}{\mathbb{R}}

\newtheorem{teo}{Theorem}[section]

\newtheorem{rema}[teo]{Remmark}
\newenvironment{dem}{\textbf Proof: }{}
\begin{document}
\title{Approximations of a complex Brownian motion with processes
constructed from a process with independent increments}
\date{}
\author{Xavier Bardina\footnote{X. Bardina is supported by the grant MTM2012-33937 from SEIDI,  Ministerio de Economia y Competividad.} and Carles Rovira\footnote{ C. Rovira is supported by the grant MTM2012-31192 from SEIDI,  Ministerio de Economia y Competividad.}}

\maketitle

$^*${\rm Departament de Matem\`atiques, Facultat de Ci\`encies, Edifici C, Universitat Aut\`onoma de Barcelona, 08193 Bellaterra}.
{\it Email: }{\tt Xavier.Bardina@uab.cat}  {\it Tel:(34) 935868563}
\newline
$\mbox{ }$\hspace{0.1cm}
$^\dagger${\rm Facultat de Matem\`atiques, Universitat de Barcelona, Gran Via 585, 08007 Barcelona}.
{\it Email: }{\tt carles.rovira@ub.edu} {\it Tel:(34) 934034740}

\begin{abstract}
In this paper, we show an approximation in law of the complex
Brownian motion by processes constructed from a stochastic process with independent increments. We give sufficient conditions for the
characteristic function  of the process with independent increments that ensure the existence of the approximation. We apply these results to L\'evy
processes.
Finally we extend this results to the $m$-dimensional
complex Brownian motion.
\end{abstract}


\section{Introduction and main result.}

The purpose of this paper is to investigate a weak approximation of
a complex Brownian motion. The most typical processes taken as
approximations to Gaussian processes are usually based on 
Donsker approximations ( the functional central limit theorem) or on
Kac-Stroock type approximations. In this paper, we will deal with
this family of approximations.

In 1974, Kac \cite{K} described the solution of the telegrapher's equation in
terms of a Poisson process. Eight years later, Stroock \cite{S} showed the weak
convergence of this solution to a Brownian motion. More precisely,
given $\{N_t, t \ge 0\}$ a standard Poisson process, the laws of the
processes $x_{\varepsilon}$
$$\{x_{\varepsilon}(t)=\varepsilon\int_0^{\frac{t}{\varepsilon^{2}}}(-1)^{N_s}ds,\quad t \in [0,T]\}$$
converge weakly towards the law of  a standard Brownian motion
 in the space of continuos functions on $[0,T]$.

This result have been extended  in order to obtain approximations of other processes as, among others: $m$-dimensional Brownian process \cite{BR1},  SPDE driven by Gaussian white noise \cite{BJQ}, fractional SDE \cite{BNRT}, multiple Wiener integrals \cite{BJT} or
complex Brownian motion \cite{Ba}.

Although all this cases are built begining from a Poisson process, a detailed study of the proofs shows that the authors use only some properties of the Poisson process that can be found in a bigger class of processes as L\'evy processes. In this paper we will show how to do this extension.

More precisely, we will deal with approximations of the complex Brownian motion built from an unique stochastic process with independent increments. Let us recall that $\{B_{t},\, t\in[0,T]\}$ is a complex Brownian motion if
it takes values on  $\mathbb C$ and its real part and its imaginary part
are two independent standard Brownian motions.

We consider the processes
\begin{equation}
\{x_{\varepsilon}^{\theta}(t)=c(\theta)\varepsilon\int_{0}^{\frac{2t}{\varepsilon^2}}e^{i\theta
X_s}ds,\quad t\in[0,T]\} ,\label{eq1}
\end{equation}
where $\{X_t, t \ge 0\}$ is a stochastic processes with independent increments and $c(\theta)$ is a constant, depending on $\theta$, that we will determine latter.
Let us recall that our approximations can be written as
$$x_{\varepsilon}^{\theta}(t)=\varepsilon c(\theta)\int_{0}^{\frac{2t}{\varepsilon^2}}\cos(\theta
X_s)ds + i\varepsilon
c(\theta)\int_{0}^{\frac{2t}{\varepsilon^2}}\sin(\theta X_s)ds.$$

Notice that when $X$ is a Poisson process it has been proved in
\cite{Ba} that for $\theta \not= 0$ and $\theta \not= \pi$ the limit
is a complex Brownian motion. When $\theta=\pi$ we obtain an
alternative version of Stroock's results since in this case
$$ e^{i\theta X_s}=(-1)^{X_s}.$$

The aim of this paper is to study the weak limits of the processes (\ref{eq1}) when $\varepsilon$ tends to zero depending on the value of $\theta$, showing that L\'evy processes can be used to
approximate a complex Brownian motion.

In section 2, we recall some basic facts on L\'evy processes and we present the classical methodology to get weak approximations of Gaussian processes.

Section 3 is devoted to give the main results of the paper.  First we give some conditions on the characteristic functions of the process $X$ that ensures the weak convergence of (\ref{eq1}) to
a complex Brownan motion. Then, we show when  the characteristic functions of L\'evy processes satisfy such conditions.

In section 4 we study the $m$-dimesional case, showing how we can obtain a $m$-dimensional complex Brownian motion from an unique L\'evy process.

Along the paper $K$ denote positive constants, not depending on
$\varepsilon$,
which may change from one expression to another one.
The real part and the imaginary part of a complex number   will be
denoted by $Re[\cdot]$ and $Im[\cdot]$.

\section{Preliminaries}

\subsection{L\'evy processes}

Set $\{X_s,\, s\geq 0\}$ a L\'evy process, that is, X has stationary and independent increments, is continuous
in probability, is c\`adl\`ag and $X_0 = 0$, and it is defined on a complete probability
space $(\Omega, \cal{F},P)$. There are many important examples of L\'evy processes: Brownian
motion, Poisson Process, jump-diffusion
processes, stable processes, subordinators, etc.

Consider $\phi_{X_t}(u)$
its characteristic function. Remember that it can be written as
$$\phi_{X_t}(u)=E\left(e^{iuX_t}\right)=e^{-t \psi_X(u)},$$
where $\psi_X(u)$ is called the L\'evy exponent of $X$.

It is well known that the L\'evy exponent can be expressed, by the L\'evy-Khinchine formula, as
\begin{equation}\label{khin}\psi_{X}(u)=-aiu+\frac12\sigma^2 u^2-\int_{\R\setminus\{0\}}(e^{iux}-1-iuxI_{|x|<1})\eta(dx),\end{equation}
where $a\in\R$, $\sigma\geq0$ and $\eta$ is a L\'evy measure, that is,
$\int_{\R\setminus\{0\}}\min\{x^2,1\}\eta(dx)<\infty$.

For notation and simplicity along the paper we will write
\begin{equation}
a(u):=Re[\psi_X(u)]=\frac12\sigma^2u^2-\int_{\R\setminus\{0\}}(\cos(ux)-1)\eta(dx),\label{mm1}\end{equation}
and
\begin{equation}
b(u)=Im[\psi_X(u)]=-au-\int_{\R\setminus\{0\}}(\sin(ux)-uxI_{|x|<1})\eta(dx).\label{mm2}\end{equation}
Notice that $a(-u)=a(u)$ and $b(-u)=-b(u)$.

We refer the reader to \cite{Sato} for more information about L\'evy processes.

\subsection{Weak approximations of the complex Brownian motion}\label{sweak}

For any $\varepsilon >0$, set $
\{x_{\varepsilon}(t), t\in[0,T]\}$ a complex stochastic process with $x_\varepsilon(0)=0$.
Consider $P_{\varepsilon}$ the image law of
$x_{\varepsilon}$ in the Banach space $\mathcal C([0,T],\mathbb
C)$ of continuous functions on
$[0,T]$.

In order to prove that
$P_\varepsilon$ converges weakly as $\varepsilon$ tends to
zero towards the law on
$\mathcal C([0,T],\mathbb C)$ of a complex Brownian motion we have to check  that the family
$P_{\varepsilon}$ is tight and that the law of
all possible weak limits of $P_\varepsilon$ is the law of two independent
standard Brownian motions.

In order to prove that the family $P_\varepsilon$ is tight, we  need to prove that the laws corresponding to the real part and the
imaginary part of the processes  $x_{\varepsilon}$ are tight.
Using the Billingsley criterium (see Theorem 12.3 of
\cite{B}) and that our processes are null on the origin, it suffices
to prove that there exists a constant $K$ such that for any
$s<t$
\begin{equation}\sup_{\varepsilon}\big(E((Re[x_\varepsilon(t)-x_\varepsilon(t)])^{4}) + E((Im[x_\varepsilon(t)-x_\varepsilon(t)])^{4})\big)\leq K(t-s)^{2}.\label{ademo1}
\end{equation}

The second part of the proof consists in the  identification of the limit law.
Let $\{P_{\varepsilon_n}\}_{n}$ be a subsequence of
$\{P_{\varepsilon}\}_{\varepsilon}$ (that we will also
denote by $\{P_{\varepsilon}\}$) weakly convergent to some
probability $P$. We want to see that the canonical process
$Y=\{Y_{t}(x)=:y(t)\}$ is a complex Brownian motion under the
probability $P$, that is, the real part and the imaginary
part of this process are two independent Brownian motions.
Using Paul L\'evy's theorem it suffices to prove that under
$P$, the real part and the imaginary part  of the canonical process are
both martingales with respect to the natural filtration, $\{\mathcal
F_{t}\}$, with quadratic variations $<Re[Y],Re[Y]>_{t}=t$,
$<Im[Y],Im[Y]>_{t}=t$ and covariation $<Re[Y],Im[Y]>_{t}=0$.

To  see that under $P$ the
real part and the imaginary part of the canonical process
$X$ are martingales with respect to its natural filtration $\{\mathcal
F_{t}\}$,
we have to prove that for any $s_{1}\leq s_{2}\leq\cdots\leq s_{n}
\leq s$ and for any bounded continuous function $\varphi:\mathbb C^{n}
\longrightarrow\mathbb R$,
$$E_{P}\big[\varphi(X_{s_{1}},...,X_{s_{n}})(Re[Y_{t}]-Re[Y_{s}])\big]=0,$$
$$E_{P}\big[\varphi(X_{s_{1}},...,X_{s_{n}})(Im[Y_{t}]-Im[Y_{s}])\big]=0.$$
Since $P_{\varepsilon}\stackrel{w}{\Rightarrow}P$,  we have that,
\begin{eqnarray*}
&& \lim_{\varepsilon\to0}E_{P_{\varepsilon}}
\big[\varphi(y(s_{1}),...,y(s_{n}))(Re[y(t)]-Re[y(s)])\big] \cr
&&\cr &=& E_{P}
\big[\varphi(y(s_{1}),...,y(s_{n}))(Re[y(t)]-Re[y(s)])\big],\cr
\end{eqnarray*}
and we get the same with the imaginary part. So, it suffices to see that
\begin{eqnarray}
& &\lim_{\varepsilon \to 0}E\big(\varphi(x_\varepsilon(s_{1}),...,x_\varepsilon(s_{n})) \big(Re[x_\varepsilon(t)] - Re[x_\varepsilon(s)]\big)\big)=0,\label{ademo2}\\
& &\lim_{\varepsilon \to 0}E\big(\varphi(x_\varepsilon(s_{1}),...,x_\varepsilon(s_{n}))\big( Im[x_\varepsilon(t)] - Im[x_\varepsilon(s)]\big)\big)=0.\label{ademo3}
\end{eqnarray}
To chek the quadratic variation, it is enough to
prove that for any $s_{1}\leq s_{2}\leq\cdots\leq s_{n} \leq s$ and
for any bounded continuous function $\varphi:\mathbb C^{n}
\longrightarrow\mathbb R$,
\begin{eqnarray}
& &\lim_{\varepsilon \to 0}E\big[\varphi(x_\varepsilon(s_{1})...,x_\varepsilon
(s_{n}))\big((Re[x_\varepsilon(t)]-Re[x_\varepsilon(s)])^2-(t-s)\big)\big]=0,\label{ademo4}\\
& &\lim_{\varepsilon \to 0}E\big[\varphi(x_\varepsilon(s_{1})...,x_\varepsilon(s_{n}))\big((Im[x_\varepsilon(t)]-Im[x_\varepsilon(s)])^2-(t-s)\big)\big]=0.\label{ademo5}
\end{eqnarray}
Finally to prove that $<Re[Y],Im[Y]>_{t}=0$, it suffices to
check that for any     $s_{1}\leq s_{2}\leq\cdots\leq s_{n} \leq s$
and for any bounded continuous function $\varphi:\mathbb C^{n}
\longrightarrow\mathbb R$,
\begin{equation}\lim_{\varepsilon \to 0}E\big[\varphi(x_\varepsilon(s_{1})...,x_\varepsilon(s_{n}))
(Re[x_\varepsilon(t)]-Re[x_\varepsilon(s)])
(Im[x_\varepsilon(t)]-Im[x_\varepsilon(s)])\big]=0\label{ademo6}.\end{equation}

\section{Approximations to a complex Brownian motion}

As we have explained, we built our approximations from a stochastic
process X with independent increments. We will deal with $X$ using
the study of its characteristic function  $\phi_X$. Let us introduce
a set of   usefull hypothesis ($H^{\theta}$) for the
characteristic function $\phi_X$ of a process $X$:
\begin{itemize}

\item[($H^\theta$1)] there exists a constant $K(\theta)$ such that
$$\varepsilon^2\int_{\frac{2s}{\varepsilon^2}}^{\frac{2t}{\varepsilon^2}}
\int_{\frac{2s}{\varepsilon^2}}^{y}\|\phi_{X_{y}-X_{x}}(\theta)\|dxdy
\le K(\theta) (t-s),$$


\item[($H^\theta$2)] there exists a constant $c(\theta)$ such that
$$ \lim_{\varepsilon  \to 0} \varepsilon^2c(\theta)^2\int_{\frac{2s}{\varepsilon^2}}^{\frac{2t}{\varepsilon^2}}\int_{\frac{2s}{\varepsilon^2}}^{x}
[\phi_{X_x-X_y}(\theta)+\phi_{X_x-X_y}(-\theta)]dydx=2(t-s),$$

\item[($H^\theta$3)]
$$ \lim_{\varepsilon \to 0} \varepsilon^2
\int_{\frac{2s}{\varepsilon^2}}^{\frac{2t}{\varepsilon^2}}\int_{\frac{2s}{\varepsilon^2}}^{y}
\|\phi_{X_y-X_x}(\theta)\|\|\phi_{X_x-X_{\frac{2s}{\varepsilon^2}}}(2\theta)\|dxdy = 0.$$

\end{itemize}

In the next Theorem we give some sufficient conditions on the characteristic function of the process  $\{X_s, s \ge 0\}$  to get the convergence to a
complex Brownian motion.

\begin{teo}\label{teogen}
Let $\{X_s,\, s\geq 0\}$ be a stochastic process with independent
increments and characteristic function $\phi_X$. Set $C_X=\{\theta,
{ \mbox{ such that }} \phi_X {\mbox{ satisfies }}
(H^{\theta})  \}$.

Define for any $\varepsilon>0$ and  $\theta\in C_X$
$$\{x_{\varepsilon}^{\theta}(t)=\varepsilon c(\theta) \int_{0}^{\frac{2t}{\varepsilon^2}}e^{i\theta
X_s}ds,\quad t\in[0,T]\}$$ where $c(\theta)$ is the constant given
by hypothesis $(H^{\theta}2)$.

Consider $P_{\varepsilon}^{\theta}$ the image law of
$x_{\varepsilon}^{\theta}$ in the Banach space $\mathcal C([0,T],\mathbb
C)$ of continuous functions on
$[0,T]$. Then, 
$P_\varepsilon^{\theta}$ converges weakly as $\varepsilon$ tends to
zero, towards the law on
$\mathcal C([0,T],\mathbb C)$ of a complex Brownian motion.
\end{teo}

\begin{dem}
We will follow the method explained in Subsection \ref{sweak}.

\bigskip

{\it Step1: Tightness.} We have to check (\ref{ademo1}), that is,  that there exists a constant $K(\theta)$ such that for any
$s<t$
$$\sup_{\varepsilon}\big(E(\varepsilon c(\theta)\int_{\frac{2s}{\varepsilon^2}}^{\frac{2t}{\varepsilon^2}}\cos(\theta
N_x)dx)^{4} + E(\varepsilon
c(\theta)\int_{\frac{2s}{\varepsilon^2}}^{\frac{2t}{\varepsilon^2}}\sin(\theta
N_x)dx)^{4}\big)\leq K(\theta)(t-s)^{2}.$$
From the properties of the complex numbers we have that
\begin{eqnarray}
&&\!\!E(\varepsilon c(\theta)\int_{\frac{2s}{\varepsilon^2}}^{\frac{2t}{\varepsilon^2}}\cos(\theta
X_x)dx)^{4} +
E(\varepsilon c(\theta)\int_{\frac{2s}{\varepsilon^2}}^{\frac{2t}{\varepsilon^2}}\sin(\theta
X_x)dx)^{4} \nonumber\cr &&\nonumber\cr
&\leq&2E\|x_\varepsilon^{\theta}(t)-x_\varepsilon^{\theta}(s)\|^4\nonumber\\
&=&2c(\theta)^4\varepsilon^4E\left(\int_{\frac{2s}{\varepsilon^2}}^{\frac{2t}{\varepsilon^2}}
e^{i\theta X_v}dv
\int_{\frac{2s}{\varepsilon^2}}^{\frac{2t}{\varepsilon^2}}
e^{-i\theta X_u}du\right)^2\nonumber\\
&=&2c(\theta)^4\varepsilon^4\int_{[\frac{2s}{\varepsilon^2},\frac{2t}{\varepsilon^2}]^4}
E\left(e^{i\theta[(X_{v_1}-X_{u_1})+(X_{v_2}-X_{u_2})]}\right)dv_1dv_2du_1du_2.\label{eqtight1}
\end{eqnarray}
Using that for $x_1<x_2<x_3<x_4$ and $\rho_i \in \{0,1\}$  for
$i=1,2,3,4$ with $\sum_{i=1}^4 \rho_i=2$ we can write
\begin{eqnarray*}
&& \!\! (-1)^{\rho_4} X_{x_4} + (-1)^{\rho_3} X_{x_3}  + (-1)^{\rho_2} X_{x_2}  + (-1)^{\rho_1} X_{x_1} \\
& = &   (-1)^{\rho_4} (X_{x_4} -X_{x_3}) +  \big((-1)^{\rho_4}  + (-1)^{\rho_3}\big) (  X_{x_3}  -  X_{x_2} ) \\
& & +\big((-1)^{\rho_4}   + (-1)^{\rho_3}  + (-1)^{\rho_2}\big) ( X_{x_2} - X_{x_1}),
\end{eqnarray*}
and the last expression (\ref{eqtight1}) can be written as the sum of 24 integrals of the type
\begin{eqnarray} 2c(\theta)^4\varepsilon^4\int_{\frac{2s}{\varepsilon^2}}^{\frac{2t}{\varepsilon^2}}
\int_{\frac{2s}{\varepsilon^2}}^{x_4}
\int_{\frac{2s}{\varepsilon^2}}^{x_3}
\int_{\frac{2s}{\varepsilon^2}}^{x_2}
&E\left(e^{i\theta[c_1(X_{x_4}-X_{x_3})+c_2(X_{x_3}-X_{x_2})+c_3(X_{x_2}-X_{x_1})]}\right) \nonumber\\
&  \times dx_1dx_2dx_3dx_4.\label{eqtight2} \end{eqnarray} where $c_1\in\{1,-1\}$, $c_2\in\{-2,0,2\}$
and
$c_3\in\{1,-1\}$.
Notice that since the process $X$ has independent
increments, we have that
\begin{eqnarray*}
&&E\left(e^{i\theta[c_1(X_{x_4}-X_{x_3})+c_2(X_{x_3}-X_{x_2})+c_3(X_{x_2}-X_{x_1})]}\right)\\
&=&E\left(e^{i\theta c_1(X_{x_4}-X_{x_3})}\right)E\left(e^{i\theta
c_2(X_{x_3}-X_{x_2})}\right)E\left(e^{i\theta
c_3(X_{x_2}-X_{x_1})]}\right),\end{eqnarray*} and we obtain,
\begin{eqnarray*}&&\|E\left(e^{i\theta[c_1(X_{x_4}-X_{x_3})+c_2(X_{x_3}-X_{x_2})+c_3(X_{x_2}-X_{x_1})]}\right)
\|\\&\leq&\|\phi_{X_{x_4}-X_{x_3}}(c_1\theta)\|\|\phi_{X_{x_2}-X_{x_1}}(c_3\theta)\|\\
&\leq&\|\phi_{X_{x_4}-X_{x_3}}(\theta)\|\|\phi_{X_{x_2}-X_{x_1}}(\theta)\|,
\end{eqnarray*}
where we have used that for any random variable $Z$,
$\|\phi_Z(-u)\|=\|\phi_Z(u)\|$.

So, each one of the  24 integrals of the type (\ref{eqtight2}) is
bounded by $$
c(\theta)^4\varepsilon^2\int_{\frac{2s}{\varepsilon^2}}^{\frac{2t}{\varepsilon^2}}
\int_{\frac{2s}{\varepsilon^2}}^{x_4}\|\phi_{X_{x_4}-X_{x_3}}(\theta)\|dx_3dx_4
\varepsilon^2\int_{\frac{2s}{\varepsilon^2}}^{\frac{2t}{\varepsilon^2}}
\int_{\frac{2s}{\varepsilon^2}}^{x_2}\|\phi_{X_{x_2}-X_{x_1}}(\theta)\|dx_1dx_2.
$$
Clearly, hypothesis ($H^\theta$1) completes the proof of this step.

\bigskip

{\it Step 2: Martingale property.} It is enough to check  (\ref{ademo2}) and  (\ref{ademo3}).  So, it suffices to see that
$$E\big(\varphi(x_\varepsilon^{\theta}(s_{1}),...,x_\varepsilon^{\theta}(s_{n}))\varepsilon c(\theta)\int_{\frac{2s}{\varepsilon^{2}}}
^{\frac{2t}{\varepsilon^2}}\cos(\theta X_x)dx\big)$$ and,
$$E\big(\varphi(x_\varepsilon^{\theta}(s_{1}),...,x_\varepsilon^{\theta}(s_{n}))\varepsilon c(\theta) \int_{\frac{2s}{\varepsilon^{2}}}
^{\frac{2t}{\varepsilon^2}}\sin(\theta X_x)dx\big)$$ converge to
zero when $\varepsilon$ tends to zero.

Thus, it is enough to prove that
$$\|E\big(\varphi(x_\varepsilon^{\theta}(s_1),...,x_\varepsilon^{\theta}(s_n))
\varepsilon c(\theta)\int_{\frac{2s}{\varepsilon^{2}}}
^{\frac{2t}{\varepsilon^2}}e^{i\theta X_x}dx\big)\|$$ converges to
zero when $\varepsilon$ tends to zero.

But applying the Schwartz inequality and using that the function
$\varphi$ is bounded it is enough to prove the convergence to zero
of
\begin{eqnarray*}
&&\|E\big(\varepsilon
c(\theta)\int_{\frac{2s}{\varepsilon^2}}^{\frac{2t}{\varepsilon^2}}e^{i\theta
X_x}dx\big)^2\|\\
&=&\|E\big(\varepsilon^2
c(\theta)^2\int_{\frac{2s}{\varepsilon^2}}^{\frac{2t}{\varepsilon^2}}\int_{\frac{2s}{\varepsilon^2}}^{y}
e^{i\theta(X_x+X_y)}dxdy\big)\|\\
&=&\|E\big(\varepsilon^2
c(\theta)^2\int_{\frac{2s}{\varepsilon^2}}^{\frac{2t}{\varepsilon^2}}\int_{\frac{2s}{\varepsilon^2}}^{y}
e^{i\theta(X_y-X_x)+2i\theta(X_x-X_{\frac{2s}{\varepsilon^2}})+2i\theta X_{\frac{2s}{\varepsilon^2}}}dxdy\big)\|\\
&=&\|E\big[e^{2i\theta X_{\frac{2s}{\varepsilon^2}}}\big]
\varepsilon^2
c(\theta)^2\int_{\frac{2s}{\varepsilon^2}}^{\frac{2t}{\varepsilon^2}}\int_{\frac{2s}{\varepsilon^2}}^{y}
\phi_{X_y-X_x}(\theta)\cdot\phi_{X_x-X_{\frac{2s}{\varepsilon^2}}}(2\theta)dxdy\|.
\end{eqnarray*}
Notice that this last expression can be bounded by
$$\varepsilon^2
c(\theta)^2\int_{\frac{2s}{\varepsilon^2}}^{\frac{2t}{\varepsilon^2}}\int_{\frac{2s}{\varepsilon^2}}^{y}
\|\phi_{X_y-X_x}(\theta)\|\|\phi_{X_x-X_{\frac{2s}{\varepsilon^2}}}(2\theta)\|dxdy$$
that from ($H^\theta$3) converges to zero when $\varepsilon$ goes to
zero.


\bigskip

{\it Step 3: Quadratic variations.} It is enough to check  (\ref{ademo4}) and  (\ref{ademo5}), that is  to
prove that for any $s_{1}\leq s_{2}\leq\cdots\leq s_{n} \leq s$ and
for any bounded continuous function $\varphi:\mathbb C^{n}
\longrightarrow\mathbb R$,
$$a_\varepsilon:=E\big[\varphi(x_\varepsilon^{\theta}(s_{1})...,x_\varepsilon^{\theta}
(s_{n}))\big((Re[x_\varepsilon^{\theta}(t)]-Re[x_\varepsilon^{\theta}(s)])^2-(t-s)\big)\big]$$
and
$$b_\varepsilon:=E\big[\varphi(x_\varepsilon^{\theta}(s_{1})...,x_\varepsilon^{\theta}(s_{n}))\big((Im[x_\varepsilon^{\theta}(t)]-Im[x_\varepsilon^{\theta}(s)])^2-(t-s)\big)\big]$$
converge to zero when $\varepsilon$ tends to zero.

In order to prove that $a_\varepsilon$ and $b_\varepsilon$ converge to zero, when $\varepsilon$
goes to zero, it is enough to show that $a_\varepsilon+b_\varepsilon$ and $a_\varepsilon-b_\varepsilon$
converge to zero. But,
\begin{eqnarray*}
&  &a_\varepsilon+b_\varepsilon\\
&=&E\big[\varphi(x_\varepsilon^{\theta}(s_{1})...,x_\varepsilon^{\theta}
(s_{n}))\big(\|x_\varepsilon^{\theta}(t)-x_\varepsilon^{\theta}(s)\|^2-2(t-s)\big)\big]\\
&=&E\big[\varphi(x_\varepsilon^{\theta}(s_{1})...,x_\varepsilon^{\theta}
(s_{n}))\big(\varepsilon^2c(\theta)^2\int_{\frac{2s}{\varepsilon^2}}^{\frac{2t}{\varepsilon^2}}\int_{\frac{2s}{\varepsilon^2}}^{\frac{2t}{\varepsilon^2}}
e^{i\theta(X_v-X_u)}dvdu-2(t-s)\big)\big]\\
&=&E\big[\varphi(x_\varepsilon^{\theta}(s_{1})...,x_\varepsilon^{\theta}
(s_{n})\big])E\big(\varepsilon^2c(\theta)^2\int_{\frac{2s}{\varepsilon^2}}^{\frac{2t}{\varepsilon^2}}\int_{\frac{2s}{\varepsilon^2}}^{\frac{2t}{\varepsilon^2}}
e^{i\theta(X_v-X_u)}dvdu
-2(t-s)\big)\\
&=&E\big[\varphi(x_\varepsilon^{\theta}(s_{1})...,x_\varepsilon^{\theta}
(s_{n})\big])\big[E\big(\varepsilon^2c(\theta)^2\int_{\frac{2s}{\varepsilon^2}}^{\frac{2t}{\varepsilon^2}}\int_{\frac{2s}{\varepsilon^2}}^{v}
e^{i\theta(X_v-X_u)}dudv\big)\\&&
+E\big(\varepsilon^2c(\theta)^2\int_{\frac{2s}{\varepsilon^2}}^{\frac{2t}{\varepsilon^2}}\int_{\frac{2s}{\varepsilon^2}}^{u}
e^{-i\theta(X_u-X_v)}dvdu\big)
-2(t-s)\big)\big]\\
&=&E\big[\varphi(x_\varepsilon^{\theta}(s_{1})...,x_\varepsilon^{\theta}
(s_{n})\big])\\&&\times\big[\varepsilon^2c(\theta)^2\int_{\frac{2s}{\varepsilon^2}}^{\frac{2t}{\varepsilon^2}}\int_{\frac{2s}{\varepsilon^2}}^{x}
[\phi_{X_x-X_y}(\theta)+\phi_{X_x-X_y}(-\theta)]dydx -2(t-s)\big].
\end{eqnarray*}
Clearly, ($H^\theta$2) yields that $\lim_{\varepsilon \to 0}
(a_\varepsilon+b_\varepsilon)=0.$

It remains to see that $a_\varepsilon-b_\varepsilon$ converges to zero. But
\begin{eqnarray}
a_\varepsilon-b_\varepsilon&=&E\big[\varphi(x_\varepsilon^{\theta}(s_{1})...,x_\varepsilon^{\theta}
(s_{n}))\big[\big(\varepsilon
c(\theta)\int_{\frac{2s}{\varepsilon^2}}^{\frac{2t}{\varepsilon^2}}\cos(\theta
X_x)dx\big)^2\nonumber\\\nonumber&&-\big(\varepsilon
c(\theta)\int_{\frac{2s}{\varepsilon^2}}^{\frac{2t}{\varepsilon^2}}\sin(\theta
X_x)dx\big)^2\big]\\\nonumber
&=&\frac12E\big[\varphi(x_\varepsilon^{\theta}(s_{1})...,x_\varepsilon^{\theta}
(s_{n}))\big[\big(\varepsilon
c(\theta)\int_{\frac{2s}{\varepsilon^2}}^{\frac{2t}{\varepsilon^2}}e^{i\theta
X_x}dx\big)^2\\&&+\big(\varepsilon
c(\theta)\int_{\frac{2s}{\varepsilon^2}}^{\frac{2t}{\varepsilon^2}}e^{-i\theta
X_x}dx\big)^2\big],\label{dostermes}
\end{eqnarray}
where in the last step we have used that
$2(\alpha^2-\beta^2)=(\alpha+\beta i)^2+(\alpha-\beta i)^2$. We will
show that this two last terms go to zero. For the first one we have
that,
\begin{eqnarray*}
&&\frac12E\big[\varphi(x_\varepsilon^{\theta}(s_{1})...,x_\varepsilon^{\theta}
(s_{n}))\big(\varepsilon
c(\theta)\int_{\frac{2s}{\varepsilon^2}}^{\frac{2t}{\varepsilon^2}}e^{i\theta
X_x}dx\big)^2\big]\\
&=&E\big[\varphi(x_\varepsilon^{\theta}(s_{1})...,x_\varepsilon^{\theta}
(s_{n}))\varepsilon^2
c(\theta)^2\int_{\frac{2s}{\varepsilon^2}}^{\frac{2t}{\varepsilon^2}}\int_{\frac{2s}{\varepsilon^2}}^{y}
e^{i\theta(X_x+X_y)}dxdy\big]\\
&=&E\big[\varphi(x_\varepsilon^{\theta}(s_{1})...,x_\varepsilon^{\theta}
(s_{n})) \\ & & \qquad\qquad \times  \varepsilon^2
c(\theta)^2\int_{\frac{2s}{\varepsilon^2}}^{\frac{2t}{\varepsilon^2}}\int_{\frac{2s}{\varepsilon^2}}^{y}
e^{i\theta(X_y-X_x)+2i\theta(X_x-X_{\frac{2s}{\varepsilon^2}})+2i\theta X_{\frac{2s}{\varepsilon^2}}}dxdy\big]\\
&=&E\big[\varphi(x_\varepsilon^{\theta}(s_{1})...,x_\varepsilon^{\theta}
(s_{n}))e^{2i\theta X_{\frac{2s}{\varepsilon^2}}}\big] \\ & & \qquad\qquad \times \varepsilon^2
c(\theta)^2\int_{\frac{2s}{\varepsilon^2}}^{\frac{2t}{\varepsilon^2}}\int_{\frac{2s}{\varepsilon^2}}^{y}
\phi_{X_y-X_x}(\theta)\cdot\phi_{X_x-X_{\frac{2s}{\varepsilon^2}}}(2\theta)dxdy.
\end{eqnarray*}
Notice that this last expression can be bounded by
$$K\varepsilon^2
c(\theta)^2\int_{\frac{2s}{\varepsilon^2}}^{\frac{2t}{\varepsilon^2}}\int_{\frac{2s}{\varepsilon^2}}^{y}
\|\phi_{X_y-X_x}(\theta)\|\|\phi_{X_x-X_{\frac{2s}{\varepsilon^2}}}(2\theta)\|dxdy$$
that from ($H^\theta$3) converges to zero when $\varepsilon$ goes to
zero. Following the same computations, and using that, in general,
for any random variable $Z$, $\|\phi_Z(-u)\|=\|\phi_Z(u)\|$ we obtain
the same bound and the convergence to zero, for the
second term of  expression (\ref{dostermes}).
\bigskip

{\it Step 4: Quadratic covariation.} It is enough to check
(\ref{ademo6}). Using that
$$\alpha\beta=\frac14i[(\alpha-\beta i)^2-(\alpha+\beta i)^2],$$ the term in the right side of (\ref{ademo6}) is equal
to
\begin{eqnarray*}
&&E\big(\varphi(x_\varepsilon^{\theta}(s_1),...,x_\varepsilon^{\theta}(s_n))(\varepsilon
c(\theta)\int_{\frac{2s}{\varepsilon^2}}^{\frac{2t}{\varepsilon^2}}
\cos(\theta X_x)dx)(\varepsilon
c(\theta\int_{\frac{2s}{\varepsilon^2}}^{\frac{2t}{\varepsilon^2}}
\sin(\theta X_x)dx) \big)\cr &=& \frac14i
E\big[\varphi(x_\varepsilon^{\theta}(s_{1})...,x_\varepsilon^{\theta}
(s_{n}))\big[\big(\varepsilon
c(\theta)\int_{\frac{2s}{\varepsilon^2}}^{\frac{2t}{\varepsilon^2}}e^{-i\theta
X_x}dx\big)^2\\&&-\big(\varepsilon
c(\theta)\int_{\frac{2s}{\varepsilon^2}}^{\frac{2t}{\varepsilon^2}}e^{i\theta
X_x}dx\big)^2\big].
\end{eqnarray*}
We have already shown in  the study of (\ref{dostermes})  that this  term goes to zero.

\hfill
$\Box$

\end{dem}

Let us state now the main result of the paper. We prove that the approximations built from a L\'evy process converge  to
a complex Brownian motion.

\begin{teo}\label{teopar}
Define for any $\varepsilon>0$
$$\{x_{\varepsilon}^{\theta}(t)=\varepsilon c(\theta) \int_{0}^{\frac{2t}{\varepsilon^2}}e^{i\theta
X_s}ds,\quad t\in[0,T]\}$$
where $\{X_s,\, s\geq 0\}$ is a L\'evy process with L\'evy exponent $\psi_X$ and 
$$c(u)=\sqrt{\frac{\| \psi_X(u) \|^2 }{2 Re[\psi_X(u)]}}.$$

Consider $P_{\varepsilon}^{\theta}$ the image law of
$x_{\varepsilon}^{\theta}$ in the Banach space $\mathcal C([0,T],\mathbb
C)$ of continuous functions on
$[0,T]$. Then, for $\theta$ such that $Re[\psi_X(\theta)]Re[\psi_X(2\theta)] \not= 0$,
$P_\varepsilon^{\theta}$ converges weakly as $\varepsilon$ tends to
zero, towards the law on
$\mathcal C([0,T],\mathbb C)$ of a complex Brownian motion.
\end{teo}

\begin{dem}
The results follows  as a particular case of Theorem \ref{teogen}.
It suffices to check that the characteristic function $\phi_X$ of
the L\'evy process $X$  satisfies ($H^{\theta}$) for any $\theta$
such that  $a(\theta)a(2\theta) \not= 0$ (recall definitions
(\ref{mm1}) and (\ref{mm2})).

\bigskip

{\it Proof of }($H^\theta$1): We can write
\begin{eqnarray*}
\varepsilon^2\int_{\frac{2s}{\varepsilon^2}}^{\frac{2t}{\varepsilon^2}}
\int_{\frac{2s}{\varepsilon^2}}^{y}\|\phi_{X_{y}-X_{x}}(\theta))\|dxdy&=&
\varepsilon^2\int_{\frac{2s}{\varepsilon^2}}^{\frac{2t}{\varepsilon^2}}
\int_{\frac{2s}{\varepsilon^2}}^{y}e^{-(y-x)a(\theta)}dxdy\\
\leq \frac{2}{a(\theta)}(t-s). \end{eqnarray*} Using that
$a(\theta)>0$ we  complete the proof of ($H^\theta$1).

%
%

\medskip

{\it Proof of }($H^\theta$2): Note first that
\begin{eqnarray*}
&&\varepsilon^2c(\theta)^2\int_{\frac{2s}{\varepsilon^2}}^{\frac{2t}{\varepsilon^2}}\int_{\frac{2s}{\varepsilon^2}}^{x}
\phi_{X_x-X_y}(\theta)dydx\\
&=&\varepsilon^2c(\theta)^2\int_{\frac{2s}{\varepsilon^2}}^{\frac{2t}{\varepsilon^2}}\int_{\frac{2s}{\varepsilon^2}}^{x}
e^{-(x-y)(a(\theta)+b(\theta)i)}dydx\\
&=&\varepsilon^2\frac{c(\theta)^2}{a(\theta)+b(\theta)i}\int_{\frac{2s}{\varepsilon^2}}^{\frac{2t}{\varepsilon^2}}
\big(1-
e^{-(x-\frac{2s}{\varepsilon^2})(a(\theta)+b(\theta)i)}\big)dx\\
&=&o(\varepsilon)+2(t-s)\frac{c(\theta)^2}{a(\theta)+b(\theta)i}.
\end{eqnarray*}
Following the same computations and taking into account that
$a(-\theta)=a(\theta)$, and that $b(-\theta)=-b(\theta)$ we obtain
that
\begin{eqnarray*}
&&\varepsilon^2c(\theta)^2\int_{\frac{2s}{\varepsilon^2}}^{\frac{2t}{\varepsilon^2}}\int_{\frac{2s}{\varepsilon^2}}^{x}
\phi_{X_x-X_y}(-\theta)dydx\\
&=&o(\varepsilon)+2(t-s)\frac{c(\theta)^2}{a(\theta)-b(\theta)i}.
\end{eqnarray*}
So
\begin{eqnarray*}
&&\varepsilon^2c(\theta)^2\int_{\frac{2s}{\varepsilon^2}}^{\frac{2t}{\varepsilon^2}}\int_{\frac{2s}{\varepsilon^2}}^{x}
[\phi_{X_x-X_y}(\theta)+\phi_{X_x-X_y}(-\theta)]dydx\\&=&o(\varepsilon)+2(t-s)\big(
\frac{c(\theta)^2}{a(\theta)+b(\theta)i}+\frac{c(\theta)^2}{a(\theta)-b(\theta)i}
\big)\\&=&o(\varepsilon)+2(t-s),
\end{eqnarray*}
and ($H^\theta$2) is clearly true.

\medskip

{\it Proof of }($H^\theta$3): Notice that
\begin{eqnarray*}
&&K\varepsilon^2
\int_{\frac{2s}{\varepsilon^2}}^{\frac{2t}{\varepsilon^2}}\int_{\frac{2s}{\varepsilon^2}}^{y}
\|\phi_{X_y-X_x}(\theta)\|\|\phi_{X_x-X_{\frac{2s}{\varepsilon^2}}}(2\theta)\|dxdy\\
&=&K\varepsilon^2
\int_{\frac{2s}{\varepsilon^2}}^{\frac{2t}{\varepsilon^2}}\int_{\frac{2s}{\varepsilon^2}}^{y}
e^{-(y-x)a(\theta)}e^{-(x-\frac{2s}{\varepsilon^2})a(2\theta)}dxdy\\
&\leq&K\varepsilon^2\frac{
1}{a(\theta)}\int_{\frac{2s}{\varepsilon^2}}^{\frac{2t}{\varepsilon^2}}
e^{-(x-\frac{2s}{\varepsilon^2})a(2\theta)}dx\\
&\leq& K\varepsilon^2\frac{ 1}{a(\theta)a(2\theta)},
\end{eqnarray*}
that converges to zero when $\varepsilon$ goes to zero.

\hfill
$\Box$

\end{dem}

\begin{rema}
Given a L\'evy process with characteristic function given by the L\'evy-Khinchine formula (\ref{khin}), the condition  $Re[\psi_X(\theta)]=0$ is equivalent to
$$\frac12\sigma^2\theta^2-\int_{\R\setminus\{0\}}(\cos(\theta x)-1)\eta(dx)=0,$$
that is, $\sigma=0$ and
$$\int_{\R\setminus\{0\}}(\cos(\theta x)-1)\eta(dx)=0.$$

So, the condition $Re[\psi_X(\theta)]Re[\psi_X(2\theta)] \not= 0$
can be written as $\sigma \not=0$ or
$$(\int_{\R\setminus\{0\}}(\cos(\theta x)-1)\eta(dx))(\int_{\R\setminus\{0\}}(\cos(2\theta x)-1)\eta(dx)) \not=0.$$
\end{rema}

\begin{rema}
When we consider $\{X_t, t \ge 0\}$ a standard Poisson process it is well-known that it is a L\'evy process with L\'evy exponent
$$\psi_X(u)=-(\cos(u)-1) -i \sin(u)$$
that corresponds to the L\'evy-Khinchine formula (\ref{khin}) with $a=0, \sigma=0$ and $\eta=\delta_{\{1\}}$.
Then the condition $Re[\psi_X(\theta)]Re[\psi_X(2\theta)] \not= 0$ yields that
$\theta\not= k \pi$  for any $k \ge 1$.

When  $\theta= (2k+1) \pi$ , we have that
\begin{equation}
x_{\varepsilon}^{\theta}(t)=c((2k+1) \pi)\varepsilon \int_{0}^{\frac{2t}{\varepsilon^2}}cos((2k+1) \pi
X_s)ds=\varepsilon \int_{0}^{\frac{2t}{\varepsilon^2}}(-1)^{X_s}ds,\label{part1}
\end{equation}
that is a real process that can not converge to a complex Brownian motion. Nevertheless part of the same proof done in Theorem \ref{teogen}  (steps 1 and 2 and study of $a_\varepsilon$, note that $b_\varepsilon=0$) works to prove tant  the processes defined by (\ref{part1}) converge weakly to a standard Brownian motion.

On the other hand, when  $\theta= 2k \pi$ , we have that
$$
x_{\varepsilon}^{\theta}(t)=c(2k \pi)\varepsilon \int_{0}^{\frac{2t}{\varepsilon^2}}cos( 2k \pi
X_s)ds=0.
$$
\end{rema}

\section{The $m$-dimensional case}

The aim of this section is to extend this result to a $m$-dimensional case for any $m\ge 1$. We will give the extensions
of Theorem \ref{teogen} and \ref{teopar}.

We define for any $\varepsilon>0$ and for any $1\leq j\leq m$
$$\left\{x_{\varepsilon}^{\theta_j}(t)=\varepsilon\int_{0}^{\frac{2t}{\varepsilon^2}}e^{i\theta_jX_s}ds,\,t\in[0,T]\right\},$$
where $\{X_s,\,s\geq0\}$ is a stochastic process with independent increments and we
consider
$$\left\{x_{\varepsilon}^{\theta}(t)=\left(x_{\varepsilon}^{\theta_1},\dots,x_{\varepsilon}^{\theta_m}\right)(t),\,t\in[0,T]\right\}.$$

In order to simplify calculus and notation we will denote by $\theta$ the $m$ values $\theta_1,
\theta_2,\dots,\theta_m$ . Since we have to control more quadratic covariations we will need to introduce new hypothesis on $\theta$,
 ($\bar H^{\theta_j,\theta_h}$) for a characteristic function $\phi_X$:

$(\bar H^{\theta_j,\theta_h})$  For any $c_1 \in  \{-1,1\}$
$$ \lim_{\varepsilon \to 0} \varepsilon^2
\int_{\frac{2s}{\varepsilon^2}}^{\frac{2t}{\varepsilon^2}}\int_{\frac{2s}{\varepsilon^2}}^{y}
\|\phi_{X_y-X_x}
(\theta_j)\|\|\phi_{X_x-X_{\frac{2s}{\varepsilon^2}}}(\theta_j +c_1
\theta_h )\|dxdy = 0.$$

Then, the extension of Theorem  \ref{teogen}, reads as follows:

\begin{teo}\label{teo2}
Let  $\{X_s,\, s\geq 0\}$ be  a
stochastic process with independent increments and characteristic
function $\phi_X$. Set $C_X^m=\{\theta \in \R^m,$  such that $\phi_X$
satisfies ($H^{\theta_j}$) for any $ j=1,\ldots,m$ and
satisfies  $ (\bar H^{\theta_j,\theta_h})$  for any $h \not= j \}$.

Define for any $\varepsilon>0$ and for any $1\leq j \leq m$
$$\{x_{\varepsilon}^{\theta_j}(t)=\varepsilon c(\theta_j) \int_{0}^{\frac{2t}{\varepsilon^2}}e^{i\theta_j
X_s}ds,\quad t\in[0,T]\},$$ 
where $c(\theta_j)$ is the constant given
by hypothesis $(H^{\theta}2)$.

Consider $P_{\varepsilon}^{\theta}$ the image law of
$x_{\varepsilon}^{\theta}=\left(x_{\varepsilon}^{\theta_1},\dots,x_{\varepsilon}^{\theta_m}\right)$ in the Banach space $\mathcal C([0,T],\mathbb
C^m)$ of continuous functions on
$[0,T]$. Then, if
$\theta\in C_X^m$,  $P_\varepsilon^{\theta}$
converges weakly as $\varepsilon$ tends to zero towards the law on
$\mathcal C([0,T],\mathbb C^{m})$ of a $m$-dimensional complex
Brownian motion.
\end{teo}

\begin{dem} The proof follows applying the computations done for the one-dimensional case combined to the method used in \cite{BR1}. We will only give some hints of the proof.

Notice that the proof of the tightnes, the martingale property of each component ans the quadratic variations can be done following exactly the proof of the one-dimensional case. So, it remains only to study all the  covariations. As it can be seen in Section 3.1 in \cite{BR1}, it
 suffices to prove that for $j\neq h$ and  for any     $s_{1}\leq
s_{2}\leq\cdots\leq s_{k} \leq s<t$ and for any bounded continuous
function $\varphi:\mathbb C^{mk} \longrightarrow\mathbb R$,

\begin{eqnarray*}
& &E\left(\varphi\left(x_{\varepsilon}^{\theta}(s_1),\dots,x_{\varepsilon}^{\theta}(s_k)\right)
\left(\varepsilon\int_{\frac{2s}{\varepsilon^2}}^{\frac{2t}{\varepsilon^2}}c(\theta_j)\cos(\theta_j
X_x)dx\right)  \right.\\  & & \qquad  \qquad \qquad \qquad  \qquad \qquad \times \left.
\left(\varepsilon\int_{\frac{2s}{\varepsilon^2}}^{\frac{2t}{\varepsilon^2}}c(\theta_h)\cos(\theta_h
X_y)dy\right)\right),
\\
& &E\left(\varphi\left(x_{\varepsilon}^{\theta}(s_1),\dots,x_{\varepsilon}^{\theta}(s_k)\right)
\left(\varepsilon\int_{\frac{2s}{\varepsilon^2}}^{\frac{2t}{\varepsilon^2}}c(\theta_j)\sin(\theta_j
X_x)dx\right)  \right.\\  & & \qquad  \qquad \qquad \qquad  \qquad \qquad \times \left.
\left(\varepsilon\int_{\frac{2s}{\varepsilon^2}}^{\frac{2t}{\varepsilon^2}}c(\theta_h)\sin(\theta_h
X_y)dy\right)\right)
\end{eqnarray*}
and
\begin{eqnarray*}
& &
E\left(\varphi\left(x_{\varepsilon}^{\theta}(s_1),\dots,x_{\varepsilon}^{\theta}(s_k)\right)
\left(\varepsilon\int_{\frac{2s}{\varepsilon^2}}^{\frac{2t}{\varepsilon^2}}c(\theta_j)\cos(\theta_j
X_x)dx\right) \right.\\  & & \qquad  \qquad \qquad \qquad  \qquad \qquad \times \left.
\left(\varepsilon\int_{\frac{2s}{\varepsilon^2}}^{\frac{2t}{\varepsilon^2}}c(\theta_h)\sin(\theta_h
X_y)dy\right)\right)
\end{eqnarray*}
converge to zero when $\varepsilon$ tends to zero.
But, using that $\cos(\theta)=\frac{e^{i\theta}+e^{-i\theta}}{2}$
and $\sin(\theta)=\frac{e^{i\theta}-e^{-i\theta}}{2i}$ ,
and the symmetry between $x$ and $y$ (interchanging the
roles of $j$ and $h$), it is enough to show that
\begin{equation}\label{otra}
\lim_{\varepsilon \to 0} \varepsilon^2 \| E\left(\varphi\left(x_{\varepsilon}^{\theta}(s_1),\dots,x_{\varepsilon}^{\theta}(s_k)\right)
\int_{\frac{2s}{\varepsilon^2}}^{\frac{2t}{\varepsilon^2}}
\int_{\frac{2s}{\varepsilon^2}}^{y}e^{i(c_1 \theta_jX_x+c_2 \theta_hX_y)}dxdy\right)\| = 0,
\end{equation}
for any $c_1, c_2 \in \{-1,1\}.$
But,
\begin{eqnarray*}
&& \| E\left(\varphi\left(x_{\varepsilon}^{\theta}(s_1),\dots,x_{\varepsilon}^{\theta}(s_k)\right)
\int_{\frac{2s}{\varepsilon^2}}^{\frac{2t}{\varepsilon^2}}
\int_{\frac{2s}{\varepsilon^2}}^{y}e^{i(c_1\theta_jX_x+c_2\theta_h X_y)}dxdy\right)\| \\
&=&\|  E\left(\varphi\left(x_{\varepsilon}^{\theta}(s_1),\dots,x_{\varepsilon}^{\theta}(s_k)\right)
\right.   \\   &&   \left. \quad \times
\int_{\frac{2s}{\varepsilon^2}}^{\frac{2t}{\varepsilon^2}}
\int_{\frac{2s}{\varepsilon^2}}^{y}e^{ ic_2\theta_h(X_y-X_x)}
e^{i(c_1\theta_j+c_2\theta_h)(X_x-X_{\frac{2s}{\varepsilon^2}})}e^{i(c_1\theta_j+c_2\theta_h)X_{\frac{2s}{\varepsilon^2}}}
dxdy\right)\|\\
&\leq&K \int_{\frac{2s}{\varepsilon^2}}^{\frac{2t}{\varepsilon^2}}
\int_{\frac{2s}{\varepsilon^2}}^{y} \|
\phi_{(X_y-X_x)}(c_2\theta_h)\|
\|\phi_{(X_x-X_{\frac{2s}{\varepsilon^2}})}(c_1\theta_j+c_2\theta_h)
\| dxdy, \\&\leq&K
\int_{\frac{2s}{\varepsilon^2}}^{\frac{2t}{\varepsilon^2}}
\int_{\frac{2s}{\varepsilon^2}}^{y} \| \phi_{(X_y-X_x)}(\theta_h)\|
\|\phi_{(X_x-X_{\frac{2s}{\varepsilon^2}})}(\theta_j+c_3\theta_h) \|
dxdy,
\end{eqnarray*}
for $c_3\in\{-1,1\}$ and (\ref{otra}) follows from  ($\bar
H^{\theta_j,\theta_h}$). \hfill $\Box$

\end{dem}

Finally, we state the extension of  Theorem \ref{teopar}.

\begin{teo}\label{te4r}
Assume now that  $\{X_s,\, s\geq 0\}$ is a L\'evy process
 with L\'evy exponent $\psi_X$ and set
$$c(u)=\sqrt{\frac{\| \psi_X(u) \|^2 }{2 Re[\psi_X(u)]}}.$$

Consider $P_{\varepsilon}^{\theta}$ the image law of
$x_{\varepsilon}^{\theta}$ in the Banach space $\mathcal
C([0,T],\mathbb C^m)$ of continuous functions on $[0,T]$. Then, for
$\theta$ such that $Re[\psi_X(\theta_j)]Re[\psi_X(2\theta_j)] \not=
0$   for all $j \in \{1,\ldots,m\}$ and $Re[\psi_X(\theta_j + c_1
\theta_h)]\not=0$ for all $j,h \in \{1,\ldots,m\}$ and $c_1 \in
\{-1,1\}$, $P_\varepsilon^{\theta}$ converges weakly as
$\varepsilon$ tends to zero, towards the law on $\mathcal
C([0,T],\mathbb C^m)$ of a $m$-dimensional complex Brownian motion.
\end{teo}

\begin{dem}
As in the one-dimensional case it suffices to check that the
characteristic function $\phi_X$ of the L\'evy process satisfies
($H^{\theta_j}$) for any $ j=1,\ldots,m$ and satisfies $(\bar
H^{\theta_j,\theta_h})$ for any $ j \not= h$. It remains only to see
the second part and can be easily checked that, for $c_1 \in \{-1,1\}$
\begin{eqnarray*}
&& \int_{\frac{2s}{\varepsilon^2}}^{\frac{2t}{\varepsilon^2}}
\int_{\frac{2s}{\varepsilon^2}}^{y} \| \phi_{(X_y-X_x)}(\theta_j) \|
\| \phi_{(X_x-X_{\frac{2s}{\varepsilon^2}})}(\theta_j+c_1 \theta_h)
\|
dxdy\\
&\leq&K\varepsilon^2\frac{1}{a(\theta_j)a(\theta_j+c_1 \theta_h)}.
\end{eqnarray*}

\hfill
$\Box$

\end{dem}

\end{document}